\documentclass[reqno]{amsart}%
\usepackage{amsfonts}
\usepackage{verbatim}
\usepackage{xcolor,soul}
\usepackage{amsmath}
\usepackage{amssymb}
\usepackage{graphicx,caption}
\usepackage{accents}
\usepackage{graphicx}%
\colorlet{cyan2}{cyan!30}
\sethlcolor{cyan2}

\providecommand{\U}[1]{\protect\rule{.1in}{.1in}}

\numberwithin{equation}{section}
\begin{document}
\title[Shallow water system]{Inverse problem for the nonlinear long wave runup on a plane sloping beach}
\author{Alexei Rybkin}
\address{Alexei Rybkin, Department of Mathematics and Statistics, University of Alaska
Fairbanks, PO Box 756660, Fairbanks, AK 99775.}
\email{arybkin@alaska.edu}
\author{Efim Pelinovsky}
\address{Efim Pelinovsky, NSE and Institute of Applied Physics, Nizhny Novgorod, Russia}
\email{Pelinovsky@gmail.com}
\author{Noah Palmer}
\address{Noah Palmer, Department of Mathematics and Statistics, University of Alaska
Fairbanks, PO Box 756660, Fairbanks, AK 99775}
\email{njpalmer@alaska.edu}
\thanks{AR and NP are supported in part by the NSF grant DMS-2009980. EP is partially
supported by Laboratory of Dynamical Systems and Applications NRU HSE, of the
Ministry of science and higher education of the RF grant ag. 075-15-2022-1101. }
\date{July, 2023}

\begin{abstract}
We put forward a simple but effective explicit method of recovering initial
data for the nonlinear shallow water system from the reading at the shoreline.
We then apply our method to the tsunami waves inverse problem.

\end{abstract}
\maketitle

\section{Introduction}

The process of the tsunami wave run-up on a coast, a very important problem of
the tsunami science, is commonly studied in the framework of shallow water
approximations \cite{LevinNosov2016}. A similar problem occurs when long swell
waves come ashore forming rogue waves \cite{Didenkulova et al 2011}. From the
mathematical point of view it is the Cauchy problem for a system of nonlinear
PDEs with initial conditions specified at the source of the tsunami wave (e.g.
epicenter of an earthquake). That is, given initial water displacement and
velocity field, find the motion of the shoreline. Unfortunately, even
displacement data are never available and certain models for initial data are
used instead (see e.g. \cite{Okada92}).

At the same time, mareographs installed at most of ports across the globe
collect data in digital form. This readily suggests an inverse problem: given
data read at specific locations, recover the characteristics of the source of
a tsunami wave. Historically this problem was approached first using isochron
\cite{Miyabe34} (see \cite{Kaistrenko72,Piatanesi01,Pires03,Voronin15} and
particularly the review by K. Satake \cite{Satake} for contemporary numerical
methods and the literature cited therein).

This note is concerned with an inverse problem for the long wave run-up on a
plane sloping beach. There is an extensive literature devoted to the direct
problem, i.e. finding various explicit solutions of the nonlinear
shallow-water equations and computing runup characteristics from the known
initial wave in a source. See e.g. \cite{Antuono07}-\cite{Dobrokhotov17}%
,\cite{Hartle21},\cite{Kanoglu04}-\cite{Madsen10},\cite{Nicolsky18}%
,\cite{Petersen2021},\cite{Rybkin14}-\cite{Synolakis87} and the literature
cited therein. In our work we solve the inverse problem when the source
characteristics can be found via the moving shoreline. More specifically,
given the time of an earthquake and the equation of motion of the shoreline
(shoreline equation), estimate the dimensions of the tsunami source. This
would not be much of a problem in the linear theory but the run-up process is
essentially nonlinear, which presents a real challenge as inverse problems in
nonlinear settings are particularly poorly understood from both geophysical
and mathematical points of view. What we do in this contribution is to find a
model that can be effectively linearized by a suitable hodograph
transformation while still retaining most important features. We take a
popular model of the 1D nonlinear shallow water system for an (infinite)
sloping beach and show that the well-known Carrier-Greenspan hodograph turns
it into a linear model whose inverse problem can be solved explicitly by means
of the Abel transform under mild additional assumption that are acceptable (in
fact, common) from a physical point of view. This provides a convenient
setting to identify and treat much more complicated models as well as a quick
assessment of what needs to be known and done in more general situations.

\section{Shallow water equations (SWE)}

\begin{figure}[h]
\centering
\captionsetup{width=\linewidth}
\includegraphics[width=\linewidth]{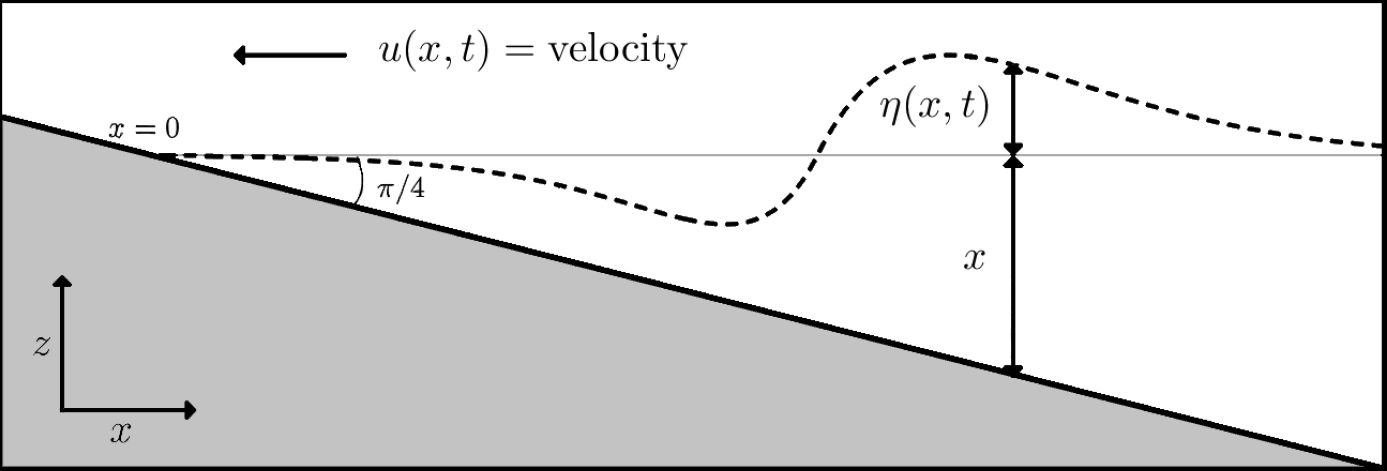}\caption{A sketch of the beach
in the $x$-$z$ plane where $\eta\left(  x,t\right)  $ is the water elevation
over the unperturbed water level $z=0$ and $H(x,t)=x+\eta(x,t)$ is the total
perturbed water depth. The unperturbed shoreline occurs at $x=0$, $z=0$.}%
\label{Fig:1}%
\end{figure}

A tsunami wave is a motion of viscous fluid described by the Navier-Stokes
equations, a highly nonlinear 3+1 (three spatial and one temporal derivatives)
system, which is notoriously hard to analyze even numerically. Under certain
assumptions (e.g. no vorticity, no friction, no dispersion, 1D velocity field,
small depth to wavelength ratio, etc.) this system can be simplified to a 2+1
system of order one equations referred to as the shallow water equations (SWE)
\cite{Johnson}. We also assume that our bathymetry represents an infinite
sloping plane beach extending along the $y$ axis infinitely far (see Figure
\ref{Fig:1}). \newline The SWE can be reduced down to a 1+1 system which in
dimensionless variables reads \cite{Johnson}%
\begin{equation}%
\begin{array}
[c]{ccc}%
\partial_{t}\eta+\partial_{x}[(x+\eta)u]=0 &  & \text{(continuity equation)}\\
\partial_{t}u+u\partial_{x}u+\partial_{x}\eta=0 &  & \text{(momentum
equation)}%
\end{array}
,\label{SWE1}%
\end{equation}
where: $\eta\left(  x,t\right)  $ is the water elevation over unperturbed
water level $z=0$ (unperturbed water depth). It need not be sign definite (can
be positive or negative); $u\left(  x,t\right)  $ is the flow velocity
averaged over the $z$ axis. Since the positive $x$-axis is directed off-shore,
$u<0$ and $u>0$ typically\footnote{In some cases an incoming wave may generate
regions where $u>0.$} corresponds respectively to an in-coming wave (i.e.
moving towards the shore) and out-going wave (i.e. moving from the shore).

Note that in the physical literature (\ref{SWE1}) is typically given in
dimensional units with the acceleration due to gravity $g$ and the slope
$\alpha$ are present explicitly. The substitution $\widetilde{x}=\left(
H_{0}/\alpha\right)  \ x,$ $\ \widetilde{t}=\sqrt{H_{0}/g}\ t/\alpha,$
$\widetilde{\eta}=H_{0}\eta,$ $\ \widetilde{u}=\sqrt{H_{0}g}\ u,$where $H_{0}$
is a typical height (or length), transforms our system (\ref{SWE1}) into a
dimensional one.

Since $H(x,t)=x+\eta\left(  x,t\right)  $ is the total (perturbed) water depth
along the main axis $x$. From Figure $1$, the equation%
\begin{equation}
x+\eta\left(  x,t\right)  =0\text{ \ (shoreline equation in physical plane)}
\label{shoreline}%
\end{equation}
describes the motion of the shoreline (wet/dry boundary) and its solution
$x_{0}\left(  t\right)  $ describes the run-up and run-down (draw down) of
tsunami waves.

We are concerned with the initial value problem (IVP) for (\ref{SWE1}) with%
\begin{equation}%
\begin{array}
[c]{ccc}%
\eta\left(  x,0\right)  =\eta_{0}\left(  x\right)  ,\ \ u\left(  x,0\right)
=u_{0}\left(  x\right)   &  & \text{(initial conditions)}%
\end{array}
\label{IC1}%
\end{equation}
and refer to such initial conditions (IC) as standard.

The IVP (\ref{SWE1})-(\ref{IC1}) has drawn an enormous attention (see, the
literature cited in Introduction). The system (\ref{SWE1}) has a quadratic
nonlinearity but its main feature is that it can be linearized by using the
Carrier-Greenspan (CG) hodograph (transformation), originally introduced in
\cite{Carrier58}. We use its version given in \cite{Carrier03}:%
\begin{equation}%
\begin{array}
[c]{cc}%
\varphi\left(  \sigma,\tau\right)  =u\left(  x,t\right)  , & \psi\left(
\sigma,\tau\right)  =\eta\left(  x,t\right)  +u^{2}\left(  x,t\right)  /2\\
\sigma^{2}=x+\eta\left(  x,t\right)  , & \tau=\left(  t-u\left(  x,t\right)
\right)  /2
\end{array}
, \label{CG}%
\end{equation}
which turns (\ref{SWE1}) into a linear hyperbolic system which, in turn, can
be reduced to the wave equation%
\begin{equation}
\partial_{\tau}^{2}\psi=\partial_{\sigma}^{2}\psi+\sigma^{-1}\partial_{\sigma
}\psi. \label{wave eq}%
\end{equation}
Note that the simple scaling $\tau=\lambda/2$ turns (\ref{wave eq}) into the
one from \cite{Carrier03}. Our form is more convenient for our purposes.

The IC (\ref{IC1}) transforms as follows. Let $x=\gamma\left(  \sigma\right)
$ solve the equation $x+\eta_{0}\left(  x\right)  =\sigma^{2}$. In the
hodograph plane $\left(  \sigma,\tau\right)  $ (\ref{IC1}) turns into%
\begin{equation}
\left.  \varphi\right\vert _{\Gamma}=\varphi_{0}(\sigma),\left.
\psi\right\vert _{\Gamma}=\psi_{0}(\sigma), \label{IC}%
\end{equation}
where $\Gamma=\left\{  \left(  \sigma,-u_{0}\left(  \gamma\left(
\sigma\right)  \right)  /2\right)  |\ \ \sigma\geq0\right\}  $ is a curve in
the $(\sigma,\tau)$ plane and $\varphi_{0}(\sigma)=u_{0}\left(  \gamma\left(
\sigma\right)  \right)  ,\ \ \ \psi_{0}(\sigma)=\eta_{0}\left(  \gamma\left(
\sigma\right)  \right)  +u_{0}\left(  \gamma\left(  \sigma\right)  \right)
^{2}/2.$

Thus, the nonlinear system (\ref{SWE1}) with linear IC (\ref{IC1}) is
transformed into a linear wave equation (\ref{wave eq}) initialized on a
(nonlinear) curve \cite{Rybkin19,Rybkin21}; $\sigma^{2}$ being the wave height
from the bottom, $\tau$ is a delayed time, $\varphi$ being the flow velocity,
and $\psi$ being wave energy. As is well-known, the main advantage of the CG
hodograph is that moving shoreline (wet/dry boundary) $x_{0}\left(  t\right)
$ given by (\ref{shoreline}) is fixed now at $\sigma=0$; the main drawback is
that the curve is nonlinear and common methods for solving IVP for linear
systems fail. The latter is not an issue when $u_{0}\left(  x\right)  =0$ (no
initial velocity). Indeed, $\Gamma$ becomes a vertical line $\left(
\sigma,0\right)  $ and the IC in (\ref{IC}) read%
\begin{equation}
\varphi(\sigma,0)=0,\ \ \ \psi(\sigma,0)=\eta_{0}\left(  \gamma\left(
\sigma\right)  \right)  .\label{LIC}%
\end{equation}
Note that in (\ref{SWE1}) $t\geq0$ but $x$ is not sign definite while in
(\ref{wave eq}) $\sigma\geq0$ but $\tau$ need not be positive. Also note that
$\sigma=0$ is a regular singular for (\ref{wave eq}) causing some
computational difficulties at the shoreline. Finally we emphasize that the CG
hodograph works as long as it is invertible, i.e. when $\det\dfrac
{\partial\left(  \sigma,\tau\right)  }{\partial\left(  x,t\right)  }\neq0$.
Recall (see e.g. \cite{Rybkin21}) that the latter holds as long as the wave
does not break. We shall always assume this condition.

\section{The shoreline equation}

In this section we are concerned with what we call the shoreline equation, an
explicit equation relating $\varphi\left(  0,\tau\right)  $, $\psi\left(
0,\tau\right)  $ with $\varphi\left(  \sigma,0\right)  =\varphi_{0}\left(
\sigma\right)  $, $\psi\left(  \sigma,0\right)  =\psi_{0}\left(
\sigma\right)  $. This is of course a classical Cauchy problem but $\tau$ need
not be positive. Using standard techniques of the Hankel transform we can
derive%
\begin{equation}
\psi(0,\tau)=\frac{d}{d\tau}\operatorname*{sgn}\tau\left[  \int_{0}%
^{\lvert\tau\rvert}\frac{s\psi_{0}\left(  s\right)  }{\sqrt{\tau^{2}-s^{2}}%
}dx-\tau\int_{0}^{\lvert\tau\rvert}\frac{s\varphi_{0}(s)}{\sqrt{\tau^{2}%
-s^{2}}}ds\right]  , \label{shoreline eq}%
\end{equation}
which is the classical Poisson formula adjusted to the signed time $\tau$.
These formulas can be conveniently written in terms of the Abel transform pair
\cite{Deans}
\[
(\mathcal{A}f)(x)=%
{\displaystyle\int_{0}^{x}}
\frac{f(s)ds}{\sqrt{x^{2}-s^{2}}},\ \ \ (\mathcal{A}^{-1}f)(x)=\frac{2}{\pi
}\dfrac{d}{dx}\int_{0}^{x}\frac{sf(s)ds}{\sqrt{x^{2}-s^{2}}}.
\]
Since negative $\tau$ in (\ref{shoreline eq}) may occur only for over critical
flows \cite{Antuono07,Antuono10} (e.g. when an extreme run-down at the initial
moment occurs) we can safely assume that $\tau\geq0$. Then%

\begin{align*}
\psi(0,\tau) &  =\frac{d}{d\tau}%
{\displaystyle\int_{0}^{\tau}}
\frac{s\psi_{0}\left(  s\right)  }{\sqrt{\tau^{2}-s^{2}}}dx-\tau\frac{d}%
{d\tau}%
{\displaystyle\int_{0}^{\tau}}
\frac{s\varphi_{0}(s)}{\sqrt{\tau^{2}-s^{2}}}ds-\int_{0}^{\tau}\frac
{s\varphi_{0}(s)}{\sqrt{\tau^{2}-s^{2}}}ds\\
&  =\frac{\pi}{2}\left[  \mathcal{A}^{-1}\psi_{0}-\tau\mathcal{A}^{-1}%
\varphi_{0}\right]  -\mathcal{A}\left(  s\varphi_{0}\right)  .
\end{align*}
Applying the Abel transform, we have%
\[
\left(\mathcal{A}\Psi\right)(\sigma)\mathcal{=}\frac{\pi}{2}\left[  \psi_{0}(\sigma)-\mathcal{A}%
\tau\left(\mathcal{A}^{-1}\varphi_{0}\right)(\sigma)\right]  -\left(\mathcal{A}^{2}  s\varphi_{0}\right)(\sigma) ,
\]
where we recall that $\psi_0$ is a function of $\sigma$.
We now compute the repeated Abel transforms seen above. Denote
\[
K_{n}\left(  x,s\right)  =\int_{s}^{x}\frac{t^{n}dt}{\sqrt{\left(  t^{2}%
-s^{2}\right)  \left(  x^{2}-t^{2}\right)  }},\ \ \ n=0,2,
\]
which is an elliptic integral. By switching the order of integration, one obtains
\[
\left(\mathcal{A}^{2} s\varphi_0\right)(\sigma){=}\int_0^\sigma s\varphi_0(s) K_0(\sigma,s)ds .
\]
The other repeated transformation is computed by integrating the innermost integral by parts and then switching the order of integration to obtain
\[
\mathcal{A}\tau\left(\mathcal{A}^{-1}\varphi_0\right)(\sigma){=}\sigma\varphi_0(0)+\int_0^\sigma K_2(\sigma,s)\varphi^\prime_0(s)ds .
\]
Thus
\begin{align}
&  \psi_{0}\left(  \sigma\right)  +%
{\displaystyle\int_{0}^{\sigma}}
\left[  \partial_{s}K_{2}\left(  \sigma,s\right)  -\frac{2s}{\pi}K_{0}\left(
\sigma,s\right)  \right]  \varphi_{0}\left(  s\right)  ds\label{int eq}\\
&  =\frac{2}{\pi}(\mathcal{A}\Psi)(\sigma),\text{ \ (shoreline equation).}%
\nonumber
\end{align}
This equation is fundamental to solving the inverse
problem for tsunami waves. In different forms and for specific $\psi
_{0},\varphi_{0}$ it appears in e.g. \cite{Carrier58}, \cite{Carrier03},
\cite{Kanoglu04}, \cite{Rybkin21}, \cite{Shimozono20}.

\section{Inverse problem}

We are now concerned with the following inverse problem. Given a sloping plane
bathymetry and assuming zero initial velocity\footnote{Such an assumption is
typical in the tsunami problem, when the bottom displacement is specified (in
the framework of shallow water, it is equivalent to the displacement of the
water surface) by the Okada seismic model} \cite{LevinNosov2016}.
$u_{0}\left(  x\right)  =0$, suppose we know the law of motion $x_{0}\left(
t\right)  $ of the shoreline. Find the initial displacement $\eta_{0}\left(
x\right)  $. In this section we give a complete solution to this problem. I.e.
we show that (\ref{CG}) and (\ref{int eq}) yield an explicit and totally
elementary procedure to restore $\eta_{0}\left(  x\right)  $ from
$x_{0}\left(  t\right)  $.

It follows from (\ref{CG}) that if $t=0$ then $\tau=0$ and $\varphi_{0}\left(
\sigma\right)  =\varphi\left(  \sigma,0\right)  =0$. Therefore (\ref{int eq})
simplifies to%
\begin{equation}
\psi_{0}\left(  \sigma\right)  =\left(  2/\pi\right)  (\mathcal{A}\Psi
)(\sigma). \label{abel}%
\end{equation}
Note also that $dx_{0}\left(  t\right)  /dt=v_{0}\left(  t\right)  $ is the
shoreline velocity. Apparently $v_{0}\left(  t\right)  =u\left(  x_{0}\left(
t\right)  ,t\right)  $ and at the shoreline ($\sigma=0$) the CG hodograph
(\ref{CG}) then reads%
\begin{equation}%
\begin{array}
[c]{cc}%
\varphi\left(  0,\tau\right)  =v_{0}\left(  t\right)  , & \psi\left(
0,\tau\right)  =\eta\left(  x_{0}\left(  t\right)  ,t\right)  +v_{0}\left(
t\right)  ^{2}/2\\
x_{0}\left(  t\right)  +\eta\left(  x_{0}\left(  t\right)  ,t\right)  =0, &
\tau=\left(  t-v_{0}\left(  t\right)  \right)  /2
\end{array}
. \label{CG1}%
\end{equation}
The inverse problem is now solved as follows:

\begin{enumerate}
\item Compute $\Psi\left(  \tau\right)  $ for (\ref{abel}) from the shoreline
data $x_{0}\left(  t\right)  $. It follows from (\ref{CG1}) that $\Psi\left(
\tau\right)  =\psi\left(  0,\tau\right)  $ can be (implicitly) found from%
\[
\Psi\left(  \tau\right)  =-x_{0}\left(  t\right)  +v_{0}\left(  t\right)
^{2}/2,\ \ \ \tau=\left(  t-v_{0}\left(  t\right)  \right)  /2.
\]

\item By (\ref{abel}) we find $\psi\left(  \sigma,0\right)  =\psi_{0}\left(
\sigma\right)  .$

\item Setting in (\ref{CG}) $\tau=0$ ($\sigma\geq0$) we find the (implicit)
equation for the initial displacement $\eta_{0}\left(  x\right)  $:%
\[
\eta_{0}\left(  x\right)  =\psi_{0}\left(  \sigma\right)  ,\ \ \ \sigma
=\sqrt{x+\eta_{0}\left(  x\right)  }.
\]

\end{enumerate}

The space limitation does not allow us to present numerical experiments here.
However, since the direct problem is well-verified already
\cite{Bueleretal2022} and the inverse problem does not contain additional
approximations, the latter is as accurate as the former and a numerical test
of our algorithm is unnecessary. Such numerics are nevertheless of interest
for applications to geophysics (tsunami problems) and will be given elsewhere.

\section{Conclusions}

Steps 1-3 above provide a smoothly working algorithm for solving the inverse
problem for tsunami waves under the assumption that the initial velocity is
zero. A similar algorithm should be in order for zero initial displacement and
nonzero initial velocity. The general case of when $u_{0}\left(  x\right)  $
and $\eta_{0}\left(  x\right)  $ are both nonzero\footnote{This situation
corresponds to a wave approaching the shore and is often implemented in
laboratory experiments.} is of course of great interest. However, on a sloping
bottom there is no rigorous description of the traveling wave (as opposed to a
flat bottom), and only approximate results can be obtained. We believe that
the important case $u_{0}\left(  x\right)  =-\eta_{0}\left(  x\right)
/\sqrt{x}$ and $u_{0}\left(  x\right)  ^{2}\ll1$ may be successfully treated
but the shoreline equation (\ref{int eq}) becomes a Fredholm integral equation
which no longer has an explicit solution. We hope to return to it elsewhere.

More complicated U and V shaped bathymetries may also be treated. In fact, we
believe that the inverse problem can be solved in closed form in the case when
the underlying direct problem can be linearized (see e.g. \cite{Bjornestad et
al 2017},\cite{Didenkulova Pelinovsky 2011},\cite{Rybkin14},\cite{Rybkin21}
for such cases). Considering general power shaped bays is work in progress.

The most practical problem is when the beach is plane for $x\leq L$ with some
$L>0$\footnote{One of the reasons is that we may neglect dispersion only when
the wave is close to the shoreline and catastrophic events are about to arise.
}. The inverse problem then requires to find the wave at $L$. The underlying
problem is boundary and the techniques of \cite{Antuono07},\cite{Rybkin21}
should be used to derive the shoreline equation in this case.

\end{document}